\documentclass[reqno,12pt,pdftex]{amsart}
\usepackage[top=25truemm,bottom=25truemm,left=25truemm,right=25truemm]{geometry}
\usepackage{amsmath,amsthm,amssymb,amsfonts,amscd,mathtools,tikz}
\usepackage{color}
\usepackage[all]{xy} 
\usepackage[colorlinks, linkcolor=black, citecolor=blue, urlcolor=black]{hyperref}
\numberwithin{equation}{section}
\theoremstyle{plain} 
	\newtheorem{thm}{Theorem}[section]
	\newtheorem{cor}[thm]{Corollary}
	\newtheorem{lem}[thm]{Lemma}
	\newtheorem{prop}[thm]{Proposition}
	\newtheorem{conj}[thm]{Conjecture}
\theoremstyle{definition}
	\newtheorem{defn}[thm]{Definition}

\theoremstyle{remark}
	\newtheorem{rem}[thm]{Remark}
	\newtheorem*{pf}{Proof}

\def\CC{\mathbb{C}}

\def\HH{\mathbb{H}}

\def\PP{\mathbb{P}}

\def\RR{\mathbb{R}}

\def\ZZ{\mathbb{Z}}
\def\A{\mathcal{A}}

\def\C{\mathcal{C}}
\def\D{\mathcal{D}}
\def\E{\mathcal{E}}

\def\P{\mathcal{P}}

\def\S{\mathcal{S}}

\def\ex#1{\langle #1 \rangle_{\rm ex}}
\def\cl#1{\overline{C}\left( #1 \right)}
\DeclareMathOperator{\Aut}{Aut}
\DeclareMathOperator{\End}{End}
\DeclareMathOperator{\EG}{EG}

\DeclareMathOperator{\Hom}{Hom}
\DeclareMathOperator{\Ind}{Ind}

\DeclareMathOperator{\Silt}{Silt}
\DeclareMathOperator{\Sim}{Sim}
\DeclareMathOperator{\SMC}{SMC}

\DeclareMathOperator{\Stab}{Stab}

\DeclareMathOperator{\fmod}{mod}
\DeclareMathOperator{\gldim}{gldim}
\DeclareMathOperator{\per}{per}

\begin{document}
\title[Stability Conditions and Algebraic Hearts for Acyclic Quivers]{Stability Conditions and Algebraic Hearts \\ for Acyclic Quivers}
\date{\today}
\author{Takumi Otani}
\address{Yau Mathematical Sciences Center, Tsinghua University, Haidian District, Beijing, China}
\email{otani-takumi-ta@alumni.osaka-u.ac.jp}
\email{takumi@tsinghua.edu.cn}
\author{Dongjian Wu}
\address{Department of Mathematics, Graduate School of Science, The University of Osaka, Toyonaka Osaka, Japan}
\email{wu.dongjian.7cx@osaka-u.ac.jp}
\maketitle
\begin{abstract}
We study stability conditions on the derived category of a finite connected acyclic quiver.
We prove that, for any stability condition on the derived category, its heart can be obtained from an algebraic heart by a rotation of phases.
Consequently, we establish the connectedness of the space of stability conditions.
Furthermore, we prove that every stability condition $\sigma$ admits a full $\sigma$-exceptional collection. 
\end{abstract}

\section{Introduction}
The space $\Stab(\D)$ of stability conditions on a triangulated category $\D$, introduced by Bridgeland \cite{Bri}, is an important homological invariant and possesses a wealth of structures.
Bridgeland proved that $\Stab(\D)$ carries a natural topological structure and, moreover, a complex manifold structure.
It is further expected that $\Stab(\D)$ admits a natural Frobenius structure in certain settings (cf.~\cite{Bri, BQS, HKK, IQ, Taka}).
This expectation is motivated by mirror symmetry, which is often understood as a correspondence among algebraic geometry, symplectic geometry, and the representation theory of algebras.
In order to approach this problem from the viewpoint of the representation theory of quivers, it is important to clarify the connection between a root system and the space of stability conditions. 
Moreover, a root system is closely related to algebraic hearts and to full exceptional collections in the derived category.
In this paper, we study stability conditions on the derived category of a finite connected acyclic quiver and their relation to these structures.
A stability condition on a triangulated category $\D$ consists of a group homomorphism $Z \colon K_0(\D) \longrightarrow \CC$ and a family of subcategories $\P = \{ \P(\phi) \}_{\phi \in \RR}$, which is an $\RR$-refinement of $t$-structures.
There are several ways to construct a stability condition.
If a heart $\A$ is {\em algebraic}, namely a length category with finitely many simple objects, then one can construct a stability condition so that $\A = \P(\phi, \phi + 1]$ for some $\phi \in \RR$.
It is natural to ask the converse question: Which stability conditions are obtained from an algebraic heart?
To answer this question, we consider the {\em support property} for a stability condition.
This notion was introduced by Kontsevich--Soibelman \cite{KS} in the study of stability conditions on derived Fukaya categories.
This condition provides various useful results for central charges.
Under the support property, we give a criterion for the heart of a stability condition to be algebraic (Proposition \ref{prop : criterion of algebraicity}, cf.~\cite{QW, Take}).
Based on the criterion, we study the existence of an algebraic heart for a stability condition on the derived category of an acyclic quiver.
To be more precise, let $Q$ be a finite connected acyclic quiver and $\D^b(Q)$ the derived category of finitely generated $k Q$-modules.
One can associate a root system with the acyclic quiver $Q$, and the classes of indecomposable objects in $\D^b(Q)$ correspond to real and imaginary roots by Kac's theorem.
The following theorem is the main result in this paper.
\begin{thm}[{Theorem \ref{thm : main 1}}]\label{thm : intro 1}
For any stability condition $\sigma = (Z, \P)$ on $\D^b(Q)$, there exists a real number $\theta \in \RR$ such that $\P(\theta, \theta + 1]$ is an algebraic heart.
\end{thm}
In the proof of Theorem \ref{thm : intro 1}, we analyze a cone $A(\sigma) \subset K_0(\D^b(Q)) \otimes_\ZZ \RR$ defined by using the support property (see Definition \ref{defn : fatten semistable imaginary cone}) instead of the set of semistable indecomposable objects whose image is an imaginary root.
We show that the cone $A(\sigma)$ is the union of finitely many connected closed cones, and the image under the central charge of the set $A(\sigma)$ is not dense.
This fact enables us to apply the criterion to our setting.
\bigskip
In the case of a Dynkin quiver $\vec{\Delta}$, it is known by \cite{KV, Qiu1} that any heart in $\D^b(\vec{\Delta})$ is algebraic and can be obtained by iteration of simple tilts from the standard heart $\fmod(k \vec{\Delta})$.
Therefore, it follows that $\Stab(\D^b(\vec{\Delta}))$ is connected.
For a finite connected acyclic quiver $Q$, it was essentially proved by Aihara--Iyama that the algebraic exchange graph of $\D^b(Q)$ is connected (See \cite{AI} and Proposition \ref{prop : connectedness of the algebraic exchange graph}).
Combining this fact with Theorem \ref{thm : intro 1}, it follows that the heart of any stability condition can be obtained by rotation and iteration of simple tilts from the standard heart $\fmod(k Q)$.
As a consequence, we obtain the following:
\begin{thm}[{Theorem \ref{thm : simple tilts and rotation}}]
$\Stab(\D^b(Q))$ is connected.
\end{thm}
It is known by Macr\`{i} that the extension closure of a full Ext-exceptional collection forms an algebraic heart.
Motivated by his work, Dimitrov--Katzarkov introduced a notion of a full $\sigma$-exceptional collection to investigate the topological structure of the space of stability conditions \cite{DK1, DK2, DK3}.
A fundamental problem in the study of full $\sigma$-exceptional collections is to establish existence.
There are several cases in which existence is known for all stability conditions: 
\begin{itemize}
\item The affine $A_1$ quiver $A^{(1)}_{1,1}$ (equivalently, the $2$-Kronecker quiver $K_2$) by \cite{Oka,Mac}.
\item The generalized Kronecker quivers $K_\ell$ with $\ell \ge 3$ by \cite{Mac, DK1}.
\item The affine $A_2$ quiver $A^{(1)}_{1,2}$ by \cite{DK1,RW}.
\item The Dynkin quivers $\vec{\Delta}$ by \cite{Ota}.
\end{itemize}
In \cite[Conjecture 3.11]{Ota}, the same statement was conjectured for affine Dynkin quivers.
Based on the correspondence between algebraic hearts and full Ext-exceptional collections, we obtain the following:
\begin{thm}[{Theorem \ref{thm : stability condition and full exceptional collections}}]\label{thm : intro 3}
Every stability condition $\sigma$ on $\D^b(Q)$ admits a monochromatic full $\sigma$-exceptional collection.
\end{thm}
This theorem not only gives an affirmative answer to the conjecture but also generalizes all known results. 
We hope that Theorem \ref{thm : intro 3} will play an important role in the study of the topological and complex structures on the space of stability conditions.
\bigskip
We briefly outline the contents of the paper. 
In Section \ref{sec : stability condition}, we recall basic definitions and properties of stability conditions on a triangulated category.
Section \ref{sec : root system} reviews root systems associated with acyclic quivers and Kac's theorem.
Section \ref{sec : main} contains the main results.
We first explain a criterion for the heart of a stability condition to be algebraic (Proposition \ref{prop : criterion of algebraicity}).
Next, we state our first main theorem (Theorem \ref{thm : main 1}), which is proved in the next subsection.
Finally, we show the connectedness of the space of stability conditions (Theorem \ref{thm : simple tilts and rotation}) and the existence of a full $\sigma$-exceptional collection (Theorem \ref{thm : stability condition and full exceptional collections}).
\bigskip
\noindent
{\bf Acknowledgements.}
The first-named author would like to thank Osamu Iyama for helpful discussions on the correspondence between silting objects and algebraic hearts.
The authors are grateful to Yu Qiu for valuable comments and suggestions.
We also thank Fabian Haiden and Atsushi Takahashi for their comments.
T.O. is supported by Beijing Natural Science Foundation Grant number IS24008.
D.W. is supported by JSPS KAKENHI KIBAN(S) 21H04994.

\section{Stability condition}\label{sec : stability condition}
Following \cite{Bri}, we recall basic notions and results for stability conditions on a triangulated category in the section.
Let $k$ be an algebraically closed field.
Throughout this paper, we always assume that our triangulated categories are $k$-linear and of finite type.

\subsection{Stability condition}
Let $\D$ be a triangulated category.
Denote by $K_0(\D)$ the Grothendieck group of $\D$.
For a full subcategory $\S \subset \D$, the extension closure is denoted by $\ex{\S}$.
\begin{defn}[{\cite[Definition 1.1]{Bri}}]
A {\em stability condition} $(Z, \P)$ on $\D$ consists of a group homomorphism $Z \colon K_0(\D)\longrightarrow\CC$ called the {\em central charge}, and a family of full additive subcategories $\P = \{ \P(\phi) \}_{\phi \in \RR}$, called the {\em slicing}, satisfying the following axioms:
\begin{enumerate}
\item For a nonzero object $E \in \P(\phi)$, we have $Z(E) = m(E) \exp(\sqrt{-1} \pi \phi)$ for some $m(E) \in \RR_{>0}$.
\item We have $\P (\phi + 1) = \P(\phi) [1]$ for all $\phi \in \RR$.
\item If $\phi_1 > \phi_2$ and $A_i \in \P(\phi_i)$ then $\Hom_\D (A_1, A_2) = 0$.
\item For each nonzero object $E \in \D$ there exists a finite sequence of real numbers
\begin{equation*}
\phi_1 > \phi_2 > \dots > \phi_n
\end{equation*}
and a collection of triangles
\begin{equation*}
\xymatrix{
0 = E_0 \ar[rr] & & E_1 \ar[r] \ar[ld] & \cdots \ar[r] & E_{n-1} \ar[rr] & & E_n = E \ar[ld] \\
& A_1 \ar@{-->}[ul] & & & & A_n \ar@{-->}[ul] & 
}
\end{equation*}
with nonzero object $A_i \in \P(\phi_i)$ for all $i = 1, \dots, n$.
\end{enumerate}
\end{defn}
The nonzero objects of $\P(\phi)$ are said to be {\em $\sigma$-semistable} of phase $\phi$, and simple objects of $\P(\phi)$ are said to be {\em $\sigma$-stable}.
Denote by $\phi(E)$ the phase of a $\sigma$-semistable object $E \in \D$.
For an object $E \in \D$ with the Harder--Narasimhan filtration as in the definition, the objects $( A_1, \cdots, A_n )$ are called its {\em Harder--Narasimhan factors} of $E$.
For any interval $I \subset \RR$, we put $\P(I) \coloneqq \ex{\P(\phi) \mid \phi \in I}$.
Then, the full subcategory $\P(0, 1]$ is a heart in $\D$, hence an abelian category.
We call $\P(0,1]$ the heart of the stability condition $\sigma$.
For a stability condition $\sigma = (Z, \P)$ on $\D$, define a subset $\C^\mathrm{ss}(\sigma) \subset K_0(\D)$ by 
\begin{equation*}
\C^\mathrm{ss}(\sigma) \coloneqq \{ \alpha \in K_0(\D) \mid \text{$\alpha = [E]$ for some $\sigma$-semistable object $E \in \D$} \}.
\end{equation*}
In order to define the support property, let us fix a norm $\| \cdot \|$ on $K_0(\D) \otimes_\ZZ \RR$.
Note that the support property does not depend on the choice of the norm.
\begin{defn}[{\cite[Definition 1]{KS}}]
We say a stability condition $\sigma = (Z, \P)$ satisfies the {\em support property} if there exists a constant $\varepsilon_\sigma > 0$ such that 
\begin{equation*}
\varepsilon_\sigma \| \alpha \| < |Z(\alpha)|
\end{equation*}
for all $\alpha \in \C^\mathrm{ss}(\sigma)$.
\end{defn}
In this paper, we always assume that our stability conditions satisfy the support property.
Denote by $\Stab(\D)$ the set of stability conditions on $\D$ with the support property.
In \cite[Section 8]{Bri}, Bridgeland introduced a natural topology on the set of stability conditions induced by a metric function.
Moreover, he also showed the forgetful map from stability conditions to central charges is a local homeomorphism, which yields a complex structure on $\Stab(\D)$ \cite[Theorem 1.2]{Bri}.
The space $\Stab(\D)$ of stability conditions has two natural actions.
The first one is the $\CC$-action defined by
\begin{equation*}
s \cdot (Z, \P) = (e^{- \pi \sqrt{-1} s} \cdot Z, \P_{{\rm Re} (s)}), \quad s \in \CC, 
\end{equation*}
where $\P_{{\rm Re} (s)}(\phi) \coloneqq \P (\phi + {\rm Re} (s))$.
The other action is given by the autoequivalence group $\Aut(\D)$:
\begin{equation*}
\Phi (Z, \P) = (Z \circ \Phi^{-1}, \Phi(\P)), \quad \Phi \in \Aut(\D).
\end{equation*}
Bridgeland gave an alternative description of a stability condition as a pair of a heart of a bounded $t$-structure and a stability function on the heart.
We recall this description here.
Note that the Grothendieck group $K_0(\A)$ of a heart $\A$ is isomorphic to $K_0(\D)$.
\begin{defn}[{\cite[Definition 2.1]{Bri}}]
Let $\A$ be a heart in $\D$.
A {\em stability function} on $\A$ is a group homomorphism $Z \colon K_0(\A) \longrightarrow \CC$ such that for all nonzero object $E \in \A$ the complex number $Z(E)$ lies in the semiclosed upper half plane 
$\HH_- \coloneqq \{ r e^{\sqrt{-1} \pi \phi} \in \CC \mid r > 0, ~ 0 < \phi \le 1 \}$.
\end{defn}
Given a stability function $Z \colon K_0(\A) \longrightarrow \CC$, the {\em phase} of a nonzero object $E \in \A$ is defined to be the real number $\phi(E) \coloneqq (1 / \pi) {\rm arg} \, Z(E) \in (0, 1]$.
A nonzero object $E \in \A$ is {\em semistable} (resp. {\em stable}) if we have $\phi(A) \le \phi(E)$ (resp. $\phi(A) < \phi(E)$) for all nonzero subobjects $A \subset E$.
We say that a stability function $Z \colon K_0(\A) \longrightarrow \CC$ satisfies the {\em Harder--Narasimhan property} if each nonzero object $E \in \A$ admits a filtration 
\begin{equation*}
0 = F_0 \subset F_1 \subset \cdots \subset F_{n - 1} \subset F_n = E
\end{equation*}
such that $F_i / F_{i - 1}$ is semistable for $i = 1, \dots, n$ with $\phi(F_1/F_0) > \phi(F_2/F_1) > \dots > \phi(F_n/F_{n - 1})$.
We say that a stability function  $Z \colon K_0(\A) \longrightarrow \CC$ satisfies the {\em support property} if there exists a constant $\varepsilon_\sigma > 0$ such that we have $\varepsilon_\sigma \| E \| < | Z(E) |$ for all semistable objects $E \in \A$.
\begin{prop}[{\cite[Proposition 5.3]{Bri}}]
To give a stability condition $\sigma = (Z, \P)$ on a triangulated category $\D$ with the support property is equivalent to giving a bounded $t$-structure on $\D$ whose heart is $\A = \P(0, 1]$ and a stability function $Z$ on its heart $\A$ with the Harder-Narasimhan property and the support property.
\qed
\end{prop}
For a heart $\A$ in $\D$, denote by $U(\A)$ the subset consisting of stability conditions on $\A$:
\begin{equation*}
U(\A) = \{ (Z, \P) \in \Stab(\D) \mid \P(0, 1] = \A \}.
\end{equation*}
Note that the subset $U(\A)$ could be empty in general.

\subsection{Global dimension function and totally semistable stability condition}\label{sec : gldim}
We recall the notion of a global dimension of a stability condition, which was defined by \cite{IQ}.
\begin{defn}[{\cite[Definition 5.4]{IQ}}]
For a slicing $\P$ on $\D$, the {\em global dimension $\gldim{\P} \in \RR_{\ge 0} \cup \{ + \infty \}$ of $\P$} is defined by
\begin{equation*}
\gldim{\P} \coloneqq \sup \{ \phi_2 - \phi_1 \in \RR \mid \Hom_\D(A_1, A_2) \ne 0 ~ \text{for} ~ A_i \in \P(\phi_i) \}.
\end{equation*}
The global dimension of a stability condition $\sigma = (Z, \P)$ on $\D$ is defined to be $\gldim{\P}$ for its slicing $\P$.
\end{defn}
For a finite-dimensional $k$-algebra $\Lambda$, let $\P_\Lambda$ denote the standard slicing given by $\P_\Lambda (1) = \fmod(\Lambda)$ and $\P_\Lambda(0, 1) = \emptyset$.
Then, we have $\gldim{\Lambda} = \gldim{\P_\Lambda}$.
Hence, the global dimension for stability conditions can be regarded as a generalization of the global dimension of a finite-dimensional algebra.
It was shown that $\gldim \colon \Stab(\D) \longrightarrow \RR_{\ge 0} \cup \{ + \infty \}$ is a continuous function.
Moreover, the global dimension function is an invariant under the $\Aut(\D)$-action and $\CC$-action.
Ikeda--Qiu showed in \cite[Lemma 5.6]{IQ} that for any stability condition $\sigma = (Z, \P)$ on $\D$ with heart $\A_\phi = \P(\phi, \phi + 1]$ we have $| \gldim{\sigma} - \gldim{\A_\phi} | \le 1$.
The following lemma is a slight modification of their statement.
\begin{lem}\label{lem : gldim for hearts}
Let $\sigma = (Z, \P)$ be a stability condition on $\D$ and $\A_\phi = \P(\phi, \phi + 1]$ for any $\phi \in \RR$.
We have 
\begin{equation*}
\gldim{\A_\phi} - 1 < \gldim{\sigma}.
\end{equation*}
\end{lem}
\begin{pf}
Let $E, F$ be objects in $\A$ and $( A_1, \cdots, A_n ), ( B_1, \cdots, B_m )$ the Harder--Narasimhan factors, respectively.
Let $\Hom_\D (E, F[k]) \ne 0$.
Then, we have $\Hom_\D (A_i, B_j[k]) \ne 0$ for some $i$ and $j$.
By definition, we have 
\begin{equation*}
k - 1 < k + \phi(B_j) - \phi(A_i) \le \gldim{\sigma}.
\end{equation*}
Since $\gldim{\A_\phi}$ takes values in $\ZZ$, we obtain the statement.
\qed
\end{pf}
A stability condition $\sigma \in \Stab(\D)$ is said to be {\em totally semistable} if every indecomposable object in $\D$ is $\sigma$-semistable.
Similarly, a totally stable stability condition is defined in a natural way.
\begin{prop}[{\cite[Proposition 3.5]{Qiu2}}]\label{prop : totally semistable stability condition}
A stability condition $\sigma$ on $\D$ is totally semistable  if and only if $\gldim{\sigma} \le 1$.
\qed
\end{prop}

\subsection{Algebraic Heart}
In this subsection, we recall the notion of an algebraic heart and collect some results for stability conditions on algebraic hearts.
This plays a central role in the paper.
An abelian category $\A$ is said to be {\em algebraic} (or {\em finite}) if it is a length category with finitely many isomorphism classes of simple objects. 
We denote by $\Sim(\A)$ the set of (isomorphism classes of) simple objects in a heart $\A$.
Note that any stability function $Z \colon K_0(\A) \longrightarrow \CC$ on an algebraic heart $\A$ satisfies the Harder--Narasimhan property.
Moreover, if the rank of the Grothendieck group is finite, a stability function $Z \colon K_0(\A) \longrightarrow \CC$ on an algebraic heart satisfies the support property.
Therefore, we have the following
\begin{prop}[{\cite[Lemma 5.2]{Bri2}, \cite[Proposition B.4]{BM}}]\label{prop : algebraic heart}
Assume that $K_0(\D) \cong \ZZ^\mu$ for some $\mu \in \ZZ_{\ge 1}$.
For an algebraic heart $\A$ with simple objects $S_1, \dots, S_\mu$, the set $U(\A)$ of stability conditions on $\A$ is isomorphic to $\HH_-^\mu$ by the map $(Z, \P) \mapsto (Z(S_1), \dots, Z(S_\mu))$.
\qed
\end{prop}
Let $\A$ be a heart in $\D$ and $S \in \A$ a simple object.
Define full subcategories
\begin{equation*}
{}^\perp S \coloneqq \{ E \in \A \mid \Hom_\A(E, S) = 0 \}, \quad
S^\perp \coloneqq \{ E \in \A \mid \Hom_\A(S, E) = 0 \}.
\end{equation*}
Then, one can consider the extension closure containing $S[1]$ and ${}^\perp S$, which is denoted by $\mu^L_S(\A)$.
It is known that the extension closure $\mu^L_S(\A)$ is a new heart in $\D$. 
We call $\mu^L_S(\A)$ the {\em left tilt} of $\A$ at $S$ (or {\em forward simple tilt} of $\A$ by $S$).
Similarly, one can define the {\em right tilted} heart $\mu^R_S(\A)$ of $\A$ at $S$ (or {\em backward simple tilt} of $\A$ by $S$). 
For an algebraic heart, the relation between stability conditions and simple tilts is described as follows:
\begin{lem}[{\cite[Lemma 5.5]{Bri2}}]\label{lem : simple tilts}
Let $\A$ be an algebraic heart in $\D$ with simple objects $S_1, \dots, S_\mu$.
If a stability condition $\sigma$ lies in the boundary of ${U}(\A)$, then either $Z(S_i) \in \RR_{> 0}$ for some $i$ and a neighbourhood of $\sigma$ is contained in $U(\A) \cup U(\mu^L_{S_i}(\A))$, or $Z(S_i) \in \RR_{< 0}$ for some $i$ and a neighbourhood of $\sigma$ is contained in $U(\A) \cup U(\mu^R_{S_i}(\A))$. 

In particular, for every $i = 1, \dots, \mu$, the unions $U(\A) \cup U(\mu^L_{S_i}(\A))$ and $U(\A) \cup U(\mu^R_{S_i}(\A))$ are connected.
\qed
\end{lem}
\begin{defn}[{\cite[Definition 5.1]{KQ}}]
The {\em exchange graph} $\EG(\D)$ of a triangulated category $\D$ is the oriented graph whose vertices are all hearts in $\D$ and whose edges correspond to left tilts between them.
We also define the {\em algebraic exchange graph} $\EG^\mathrm{alg}(\D)$ as the full subgraph of $\EG(\D)$ consisting of algebraic hearts.
\end{defn}
Denote by $\Stab^\mathrm{alg}(\D)$ the subset of $\Stab(\D)$ consisting of stability conditions whose heart is algebraic:
\begin{equation*}
\Stab^\mathrm{alg}(\D) = \bigcup_{\A \in \EG^\mathrm{alg}(\D)} U(\A).
\end{equation*}

\section{Root system and Kac's theorem}\label{sec : root system}
In this section, we recall some basic facts about root systems associated with acyclic quivers and Kac's Theorem.
We refer to \cite{Kac1, Kac2} for more details.

\subsection{Root systems associated with acyclic quivers}
Let $Q = (Q_0, Q_1)$ be a finite connected acyclic quiver and $\mu \in \ZZ_{\ge 1}$ the number of vertices.
Denote by $\overline{Q}$ the underlying graph of the quiver $Q$.
The {\em generalized Cartan matrix} $A_Q = (a_{ij})$ associated with $Q$ is defined by 
\begin{equation*}
a_{ij} \coloneqq 2 \delta_{ij} - (q_{ij} + q_{ji}),
\end{equation*}
where $q_{ij}$ is the number of arrows connecting the vertices $i$ and $j$.
By definition, the generalized Cartan matrix associated with a finite connected acyclic quiver is indecomposable and symmetric.
It is known by \cite[Theorem 4.3]{Kac2} that indecomposable generalized Cartan matrices are classified into three types, which are of finite type, affine type and indefinite type.
By the above construction, Dynkin quivers and affine Dynkin quivers correspond naturally to the finite and affine types, respectively.
Following \cite[Section 1]{Kac1}, one can associate a root system to a generalized Cartan matrix.
The {\em root lattice} is a free abelian group $L = \bigoplus_{i = 1}^\mu \ZZ \alpha_i$ with generators $\alpha_1, \dots, \alpha_\mu$, called {\em simple roots}.
Denote by $\Pi = \{ \alpha_1, \dots, \alpha_\mu \}$ the set of simple roots.
We can also define a symmetric $\ZZ$-bilinear form $I \colon L \times L \longrightarrow \ZZ$ by $I(\alpha_i, \alpha_j) = a_{ij}$.
For a simple root $\alpha_i \in \Pi$, the reflection $r_i \in \Aut_\ZZ (L, I)$ is defined by 
\begin{equation*}
r_i (\lambda) \coloneqq \lambda - I (\lambda, \alpha_i) \alpha_i, \quad \lambda \in L.
\end{equation*}
The subgroup $W \coloneqq \langle r_1, \dots, r_\mu \rangle$ of $\Aut_\ZZ (L, I)$ generated by reflections is called the {\em Weyl group}.
Let $L^+ \coloneqq \sum_{i = 1}^\mu \ZZ_{\ge 0} \alpha_i$ and $L^- \coloneqq - L^+ = \sum_{i = 1}^\mu \ZZ_{\le 0} \alpha_i$. 
Define the set of {\em real root} $\Delta_\mathrm{re}$ by 
\begin{equation*}
\Delta_\mathrm{re} \coloneqq W (\Pi) = \{ w(\alpha_i) \in L \mid w \in W, ~ i = 1, \dots, \mu \}.
\end{equation*}
The set of {\em positive real roots} $\Delta^+_\mathrm{re}$ (resp., {\em negative real roots} $\Delta^-_\mathrm{re}$) is defined by $\Delta^+_\mathrm{re} \coloneqq \Delta_\mathrm{re} \cap L^+$ (resp., $\Delta^-_\mathrm{re} \coloneqq \Delta_\mathrm{re} \cap L^-$).
Then, it is known that we have $\Delta_\mathrm{re} = \Delta^+_\mathrm{re} \sqcup \Delta^-_\mathrm{re}$ (cf.~\cite{Kac1, Kac2}).
For an element $\lambda = \sum_{i = 1}^\mu n_i \alpha_i \in L$, the {\em support} of $\lambda$ is the full subgraph of $\overline{Q}$ consisting of vertices $i \in Q_0$ for which $n_i \ne 0$.
Consider a subset $K \subset L^+ \setminus \{ 0 \}$ defined by 
\begin{equation*}
K \coloneqq \{ \lambda \in L^+ \setminus \{ 0 \} \mid \text{$\lambda$ has a connected support}, ~ I(\lambda, \alpha_i) \le 0 ~ \text{for} ~ i = 1, \dots, \mu \}.
\end{equation*}
The set of {\em positive imaginary roots} $\Delta^+_\mathrm{im}$ is defined by
\begin{equation*}
\Delta^+_\mathrm{im} \coloneqq W (K) = \{ w(\lambda) \in L^+ \mid w \in W, ~ \lambda \in K \}, 
\end{equation*}
and the set of {\em negative} imaginary roots $\Delta^-_\mathrm{im}$ is given by $\Delta^-_\mathrm{im} \coloneqq - \Delta^+_\mathrm{im}$.
Define the set of {\em imaginary roots} $\Delta_\mathrm{im}$ by $\Delta_\mathrm{im} \coloneqq \Delta^+_\mathrm{im} \sqcup \Delta^-_\mathrm{im}$.
Finally, we define the set of {\em roots} $\Delta$, {\em positive roots} $\Delta^+$ and {\em negative roots} $\Delta^-$ by 
\begin{equation*}
\Delta \coloneqq \Delta_\mathrm{re} \sqcup \Delta_\mathrm{im}, \quad \Delta^\pm \coloneqq \Delta \cap L^\pm,
\end{equation*}
respectively. 

\subsection{Indecomposable objects and Kac's theorem}\label{subsec : Kac's theorem}
Let $\D^b(Q) = \D^b (\fmod(k Q))$ be the bounded derived category of finitely generated modules over the path algebra $k Q$.
Since the path algebra $k Q$ is hereditary, we have the following (cf.~\cite[Section 4]{Hap}):
\begin{equation*}
\Ind \D^b(Q) = \bigsqcup_{p \in \ZZ} \Ind \fmod(k Q) [p],
\end{equation*}
where $\Ind \A$ denotes the set of the (isomorphism classes of) indecomposable objects in an additive category $\A$.
Let $S_i$ denote the simple $k Q$-module corresponding to the vertex $i \in Q_0$.
Since the abelian category $\fmod(k Q)$ is algebraic with simple objects $S_1, \dots, S_\mu$, we have 
\begin{equation*}
K_0(\D^b(Q)) \cong K_0(\fmod(k Q)) \cong \bigoplus_{i = 1}^\mu \ZZ [S_i],
\end{equation*}
which yields a group isomorphism
\begin{equation*}
K_0(\D^b(Q)) \longrightarrow L, \quad [S_i] \mapsto \alpha_i.
\end{equation*}
For every $i, j \in Q_0$ we have 
\begin{equation*}
\chi (S_i, S_j) + \chi(S_j, S_i) = I(\alpha_i, \alpha_j),
\end{equation*}
where $\chi \colon K_0(\D^b(Q)) \longrightarrow \ZZ$ is the Euler form defined by 
\begin{equation*}
\chi(E,F) \coloneqq \sum_{p \in \ZZ} (-1)^p \dim_k \Hom_{\D^b(Q)} (E,F[p]), \quad E, F \in \D. 
\end{equation*}
Therefore, we can identify $(K_0(\D^b(Q)), \chi + \chi^T)$ with the root lattice $(L, I)$.
By abuse of notation, the set $\Delta$ of roots is regarded as a subset of $K_0(\D^b(Q))$.
\begin{prop}[{\cite[Theorem 1]{Kac1}}]\label{prop : Kac}
An object $E \in \fmod(k Q)$ is indecomposable if and only if $[E] \in \Delta^+$.
\qed
\end{prop}
Since an indecomposable object in $\D^b(Q)$ is given as a shift of an indecomposable $k Q$-module, we have the following
\begin{cor}\label{cor : Kac}
For any indecomposable object $E \in \D^b(Q)$, we have $[E] \in \Delta$.
\qed
\end{cor}

\section{Stability conditions for acyclic quiver}\label{sec : main}
In this section, we study stability conditions on the derived category of the path algebra associated with an acyclic quiver.

\subsection{Stability condition and algebraic heart}
We first show the following proposition in a general setting, which plays an important role to study the algebraicity of a stability condition.
This proposition was already proved in \cite[Lemma 3.1]{QW} and \cite[Lemma 61]{Take}.
Nevertheless, we shall give a proof to highlight the significance of the support property.
\begin{prop}\label{prop : criterion of algebraicity}
Let $\D$ be a triangulated category.
Assume that the rank of $K_0(\D)$ is finite.
For a stability condition $(Z, \P)$, the heart $\A = \P(0, 1]$ is algebraic if and only if there is a positive number $\delta > 0$ such that $\P(0, \delta) = \{ 0 \}$.
\end{prop}
\begin{pf}
Let $\A = \P(0, 1]$ be an algebraic heart with simple objects $S_1, \dots, S_\mu$.
By Proposition \ref{prop : algebraic heart}, the real number $\delta \coloneqq \min_{i = 1, \dots, \mu} \phi(S_i) > 0$ satisfies $\P(0, \delta) = \{ 0 \}$.

Conversely, we suppose $\P(0, \delta) = \{ 0 \}$ for some $\delta > 0$.
By the $\CC$-action, we may assume that $\P(0, \delta/2) = \P(1 - \delta/2, 1) = \{ 0 \}$ and $\A = \P(0, 1]$ without loss of generality.
Then, there is a constant $M_\delta > 0$ such that for any nonzero object $E \in \A$ we have 
\begin{equation*}
M_\delta \cdot \mathrm{Im} \, Z(E)  > | \mathrm{Re} \, Z(E)|.
\end{equation*}
Now assume that $\A$ is not a length category, which implies that there is an object $E \in \A$ with an infinite composition series of simple quotients $\{ S_i \}$. 
Then we have 
\begin{equation*}
\mathrm{Im} \, Z(E) = \sum_i \mathrm{Im} \, Z(S_i) < \infty.
\end{equation*}
Note that there exists constant  $C > 0$ such that $\inf_{i} \| S_i \| > C$. 
Since any simple object in $\A$ is $\sigma$-stable, it follows from the support property that
\begin{equation*}
0 < \varepsilon_\sigma 
< \dfrac{|Z(S_i)|}{\| S_i \|} 
\le \dfrac{| \mathrm{Im} \, Z(S_i) | + | \mathrm{Re} \, Z(S_i) |}{\| S_i \|} 
\le \dfrac{(M_\delta + 1) \mathrm{Im} \, Z(S_i)}{\| S_i \|}.
\end{equation*}
Hence, we have
\begin{equation*}
 \dfrac{C \varepsilon_\sigma}{M_\delta + 1} <\mathrm{Im} \, Z(S_i),
\end{equation*}
which is a contradiction.

Therefore, every object in $\A$ has finite length. 
Consequently, the classes of simple objects in $\A$ form a basis of $K_0(\D)$, and the number of isomorphism classes of simple objects in $\A$ is equal to the rank of $K_0(\D)$.
Thus, the heart $\A$ is algebraic.
\qed
\end{pf}
Proposition \ref{prop : criterion of algebraicity} enables us to study the existence of algebraic heart for a stability condition.
In what follows, we shall consider the derived category $\D^b(Q)$ of an acyclic quiver $Q$.
The following is our main theorem in this paper.
\begin{thm}\label{thm : main 1}
Let $Q$ be a finite connected acyclic quiver.
For any stability condition $\sigma = (Z, \P)$ on $\D^b(Q)$, there is a real number $\theta \in \RR$ such that $\P(\theta, \theta + 1]$ is algebraic.
In particular, we have
\begin{equation*}
\Stab(\D^b(Q)) = \CC \cdot \Stab^\mathrm{alg}(\D^b(Q)).
\end{equation*}
\end{thm}
We shall prove Theorem \ref{thm : main 1} in Section \ref{sec : proof of main theorem 1}.
Our strategy of the proof is as follows:
Based on Kac's theorem, we analyze the image under the central charge of the set of semistable roots to apply Proposition \ref{prop : criterion of algebraicity} to our setting.
In order to observe the behavior of the semistable imaginary roots, we introduce a set $A(\sigma)$ based on the support property (Definition \ref{defn : fatten semistable imaginary cone}).
The set $A(\sigma)$ is a union of finitely many connected closed cones and includes all semistable imaginary roots (Lemma \ref{lem : fatten semistable imaginary cone} and \ref{lem : finiteness of connected components}).
Then, we show that the image under the central charge of the set $A(\sigma)$ is not dense in $\CC^\ast$ (Proposition \ref{prop : image of the fatten semistable imaginary cone}).
Finally, we deduce Theorem \ref{thm : main 1} by observing the behavior of real roots (Proposition \ref{prop : existence of gaps}).

\subsection{Proof of Theorem \ref{thm : main 1}}\label{sec : proof of main theorem 1}
Fix a finite connected acyclic quiver $Q = (Q_0, Q_1)$.
As we discussed in Section \ref{subsec : Kac's theorem}, the Grothendieck group $K_0(\D^b(Q))$ equipped with the symmetrized Euler form is identified with the root lattice $(L, I)$.
For simplicity, we will denote $L_\RR \coloneqq L \otimes_\ZZ \RR$ and $L_\RR^\ast \coloneqq L_\RR \setminus \{ 0 \}$.
Fix a norm $\| \cdot \|$ on $L_\RR$.
Then, there is a natural topology on $L_\RR$.
For a non-empty subset $S \subset L^\ast_\RR$, define a cone $C(S)$ in $L_\RR^\ast$ by 
\begin{equation*}
C(S) \coloneqq \RR_{> 0} S = \{ r \cdot \alpha \in L^\ast_\RR \mid r > 0, ~ \alpha \in S \}.
\end{equation*}
Denote by $\overline{C}(S)$ the closure of the cone $C(S)$ in $L_\RR$.
It will be convenient to put $\overline{C}(S)_0 \coloneqq \overline{C}(S) \setminus \{0\}$.
We call $\cl{\Delta_\mathrm{im}}$ the {\em imaginary cone}.
We collect some properties of the imaginary cone.
\begin{lem}[{\cite[Proposition 1.4]{Kac1}, cf.~\cite[Lemma 2.5]{Ike}}]\label{lem : positive imaginary cone}
Assume that $\Delta_\mathrm{im} \ne \emptyset$. 
Then, $\cl{\Delta_\mathrm{im}^+}_0$ is a convex cone contained in $\sum_{i = 1}^\mu \RR_{> 0} \alpha_i$.
\qed
\end{lem}
\begin{lem}[{\cite[Lemm 5.8]{Kac2}}]\label{lem : limit of real rays}
The limit rays in $L_\RR$ for $C(\Delta_\mathrm{re}^+)$ lie in $\cl{\Delta_\mathrm{im}^+}$.
\qed
\end{lem}
\begin{rem}
In \cite{Ike}, the author also considered an imaginary cone.
For a connected acyclic quiver, his imaginary cone is given by $\cl{\Delta_\mathrm{im}^+}$ in our notation.
\end{rem}
From now on, we introduce subsets and cones associated with a given stability condition.
Let $\sigma = (Z, \P)$ be a stability condition on $\D$.
Note that the central charge $Z \colon L \longrightarrow \CC$ is naturally extended as an $\RR$-linear map $L_\RR \longrightarrow \CC$, which is continuous.
For simplicity, we also write $Z \colon L_\RR \longrightarrow \CC$.
Define a continuous map $f_Z \colon L_\RR^\ast \longrightarrow [0, + \infty)$ by 
\begin{equation*}
f_Z(\alpha) \coloneqq \dfrac{|Z(\alpha)|}{\| \alpha \|}, \quad \alpha \in L_\RR^\ast.
\end{equation*}
Write $\| Z \| \coloneqq \sup_{\alpha \in L_\RR^\ast} f_Z(\alpha)$, which satisfies $\| Z \| < \infty$.
Recall that the support property yields the existence of a positive number $\varepsilon_\sigma > 0$ satisfying $f_Z(\alpha) > \varepsilon_\sigma$ for any $\alpha \in \C^\mathrm{ss}(\sigma)$.
\begin{rem}[{cf.~\cite[Remark 4.3]{Ike}}]
A stability condition $\sigma = (Z, \P)$ satisfies the support property if and only if there is no sequence $\{ \alpha_k \}_{k = 1}^\infty \subset \C^\mathrm{ss}(\sigma)$ such that $\lim_{k \to \infty} f_Z(\alpha_k) = 0$.
\end{rem}
\begin{defn}\label{defn : fatten semistable imaginary cone}
For a stability condition $\sigma = (Z, \P)$ on $\D^b(Q)$, define a subset $A(\sigma) \subset L_\RR^\ast$ as the intersection between $\cl{\Delta_\mathrm{im}}_0$ and the inverse image of the interval $[\varepsilon_\sigma, \| Z \|] \subset \RR$:  
\begin{equation*}
A (\sigma) \coloneqq \cl{\Delta_\mathrm{im}}_0 \cap f_Z^{-1}[\varepsilon_\sigma, \| Z \|].
\end{equation*}
\end{defn}
Note that it follows from the continuity of $f_Z$ that the set $A(\sigma)$ is closed in $L^\ast_\RR$.
\begin{rem}
It is known that $\Delta_\mathrm{im} = \emptyset$ for Dynkin quivers and $\Delta_\mathrm{im} = \{ n \delta \in L \mid n \in \ZZ \setminus \{ 0 \} \}$ with generator $\delta$ for affine Dynkin quivers (see \cite[Theorem 5.6]{Kac2}) .
Therefore, when $Q$ is a Dynkin quiver,  the set $A(\sigma)$ is empty for any stability condition $\sigma$ on $\D^b(Q)$.
When $Q$ is an affine Dynkin quiver, the set $A(\sigma)$ is either $\RR^\ast \delta$ or empty.
\end{rem}
\begin{lem}\label{lem : fatten semistable imaginary cone}
We have $\cl{\C^\mathrm{ss} (\sigma) \cap \Delta_\mathrm{im}}_0 \subset A (\sigma)$.
\end{lem}
\begin{pf}
Note that we have 
\begin{equation*}
\cl{\C^\mathrm{ss} (\sigma) \cap \Delta_\mathrm{im}}
\subset \cl{\C^\mathrm{ss} (\sigma)} \cap \cl{\Delta_\mathrm{im}}.
\end{equation*}
Since the map $f_Z$ is $\RR_{> 0}$-invariant and continuous, it follows from the support property that $\cl{\C^\mathrm{ss} (\sigma)}_0 \subset f_Z^{-1}[\varepsilon_\sigma, \| Z \|]$.
Hence, we have the statement.
\qed
\end{pf}
\begin{lem}\label{lem : finiteness of connected components}
The set $A (\sigma)$ has finitely many connected components.
Moreover, there are finitely many connected closed subcones $A^+_1, \cdots, A^+_n$ of $\cl{\Delta_\mathrm{im}^+}$ and $A^-_1, \cdots, A^-_n$ of $\cl{\Delta_\mathrm{im}^-}$ such that 
\begin{equation*}
A(\sigma) = A^+_1 \sqcup \cdots \sqcup A^+_n \sqcup A^-_1 \sqcup \cdots \sqcup A^-_n,
\end{equation*}
and $- A^+_i = A^-_i$ for all $i = 1, \dots, n$. 
\end{lem}
\begin{pf}
Since $\Delta_\mathrm{im} = \Delta_\mathrm{im}^+ \sqcup \Delta_\mathrm{im}^-$, Lemma \ref{lem : positive imaginary cone} implies that $\cl{\Delta_\mathrm{im}}_0 = \cl{\Delta_\mathrm{im}^+}_0 \sqcup \cl{\Delta_\mathrm{im}^-}_0$.
It follows from the connectedness of our acyclic quiver $Q$ that $\cl{\Delta_\mathrm{im}^\pm}_0$ is connected.
Hence, the set $\cl{\Delta_\mathrm{im}}_0$ has two connected components. 
On the other hand, since the map $f_Z$ is continuous, the intersection of the set $f_Z^{-1} [\varepsilon_\sigma, \| Z \|]$ and the unit sphere with respect to the norm $\| \cdot \|$ is compact.
It then follows that the set $f_Z^{-1} [\varepsilon_\sigma, \| Z \|]$ has finitely many connected components.
Hence, the number of connected components of $A (\sigma)$ is finite.
\qed
\end{pf}
\begin{lem}\label{lem : intersections of cones}
Let $A (\sigma) = A^+_1 \sqcup \cdots \sqcup A^+_n \sqcup A^-_1 \sqcup \cdots \sqcup A^-_n$ as in Lemma \ref{lem : finiteness of connected components}.
For any $i, j, k = 1, \dots, n$, we have the followings:
\begin{enumerate}
\item We have $Z(A^+_i) \cap Z(A^-_i) = \emptyset$ and $Z(A^-_i) = - Z(A^+_i)$.
\item If $Z(A^+_i) \cap Z(A^+_j) \ne \emptyset$, then $A^+_i = A^+_j$.
\item If $Z(A^+_i) \cap Z(A^-_j) \ne \emptyset$ and $Z(A^-_j) \cap Z(A^+_k) \ne \emptyset$, then $A^+_i = A^+_k$.
\item If $Z(A^-_i) \cap Z(A^-_j) \ne \emptyset$, then $A^-_i = A^-_j$.
\item If $Z(A^-_i) \cap Z(A^+_j) \ne \emptyset$ and $Z(A^+_j) \cap Z(A^-_k) \ne \emptyset$, then $A^-_i = A^-_k$.
\end{enumerate}
\end{lem}
\begin{pf}
We prove the statements one by one.

(1): Since $Z$ is $\RR$-linear and $0 \not\in Z(A^\pm_i)$, the statement is obvious.

(2): By assumption, there are $a_i \in A^+_i$ and $a_j \in A^+_j$ such that $Z(a_i) = Z(a_j)$.
For all $t \in [0, 1]$ we have 
\begin{eqnarray*}
\varepsilon_\sigma \| t a_i + (1 - t) a_j \| 
& \le & \varepsilon_\sigma t \| a_i \| + \varepsilon_\sigma (1 - t) \| a_j \| \\
& < & t |Z(a_i)| + (1 - t) |Z(a_j)| \\
& = & |Z( t a_i + (1 - t) a_j)|.
\end{eqnarray*}
Then, $a_i$ and $a_j$ lie in the same connected component of $f_Z^{-1}[\varepsilon_\sigma, \| Z \|]$.
Hence, it follows from the connectedness of $\cl{\Delta_\mathrm{im}^+}_0$ that $A^+_i = A^+_j$.

(3): There are $a_i \in A^+_i$ and $a_j \in A^-_j$ such that $Z(a_i) = Z(a_j)$, and $a'_j \in A^-_j$ and $a_k \in A^+_k$ such that $Z(a'_j) = Z(a_k)$.
Then, $a_i$ and $a_j$ lie in the same connected component of $f_Z^{-1}[\varepsilon_\sigma, \| Z \|]$ as in the argument of (2).
By definition, $a_j$ and $a'_j$ also lie in the same connected component.
Since $a'_j$ and $a_k$ belong to the same connected component, the component contains $a_i$ and $a_k$.
Therefore, since $a_i, a_k \in \cl{\Delta_\mathrm{im}^+}_0$, we finally have $A^+_i = A^+_k$.

(4) and (5) are analogues of (2) and (3), respectively.
\qed
\end{pf}
We prepare a basic lemma concerning cones in the complex plane, which is used in the proof of Proposition \ref{prop : image of the fatten semistable imaginary cone}.
We say a cone $C \subset \CC^\ast$ is generated by a subset $S \subset \CC^\ast$ if $C = \{ r \cdot s \in \CC^\ast \mid r > 0, ~ s \in S \}$.
\begin{lem}\label{lem : connected closed cone in H}
A subset $C \subset \CC^\ast$ is a convex closed cone if and only if there exist real numbers $\phi^+_C, \phi^-_C \in \RR$ such that $0 \le \phi^+_C - \phi^-_C < 1$ and $C$ is the cone generated by $\{ e^{\pi \sqrt{-1} \phi} \in \CC^\ast \mid \phi^-_C \le \phi \le \phi^+_C \}$.
\end{lem}
\begin{pf}
Our proof is based on \cite[Lemma 2.10]{Ike}.
We first consider the case $C \cap \RR_{> 0} = \emptyset$.
By the compactness of $\{ z \in C \mid |z| = 1 \}$, we have the maximum phase $\phi^+_C$ of $C$ and the minimum one $\phi^-_C$ defined by
\begin{eqnarray*}
\phi^+_C & \coloneqq \max \{ \phi \in (0, 2] \mid e^{\pi \sqrt{-1} \phi} \in C \}, \\
\phi^-_C & \coloneqq \min \{ \phi \in (0, 2] \mid e^{\pi \sqrt{-1} \phi} \in C \}.
\end{eqnarray*}
The convexity of $C$ and $0 \not\in C$ imply that $0 \le \phi^+_C - \phi^-_C < 1$.
One can check easily that $C$ is the cone generated by $\{ e^{\pi \sqrt{-1} \phi} \in \CC^\ast \mid \phi^-_C \le \phi \le \phi^+_C \}$.

Next, we assume $C \cap \RR_{> 0} \ne \emptyset$.
It follows from the convexity of $C$ and $0 \not\in C$ that $C \cap \RR_{> 0} = \emptyset$.
Then, we can show the statement in the same way.

The converse statement is obvious.
\qed
\end{pf}
For an interval $I \subset \RR$, it will be convenient to denote by $C^I$ the cone generated by $\{ e^{\pi \sqrt{-1} \phi} \in \CC^\ast \mid \phi \in I \}$.
\begin{prop}\label{prop : image of the fatten semistable imaginary cone}
There exist real numbers $\theta, \theta' \in \RR$ with $\theta' > \theta$ such that 
\begin{equation*}
Z(A(\sigma)) \cap C^{(\theta, \theta')}= \emptyset.
\end{equation*}
\end{prop}
\begin{pf}
By relabeling the indices, we may assume that there is $m \in \{ 1, \dots, n \}$ such that 
\begin{equation*}
Z(A^+_i) \cap Z(A^-_{m + i}) \ne \emptyset, \quad i = 1, \dots, n - m.
\end{equation*}
For $i = 1, \dots, m$, we define sets $C_i \subset \CC^\ast$ by 
\begin{eqnarray*}
C^+_i & \coloneqq &
\begin{cases}
Z(A^+_i) \cup Z(A^-_{m + i}), & i = 1, \dots, n - m, \\
Z(A^+_i), & i = n - m + 1, \dots, m,
\end{cases} \\
C^-_i & \coloneqq & 
\begin{cases}
Z(A^-_i) \cup Z(A^+_{m + i}), & i = 1, \dots, n - m, \\
Z(A^-_i), & i = n - m + 1, \dots, m.
\end{cases}
\end{eqnarray*}
Note that $C^-_i = - C^+_i$ for each $i = 1, \dots, m$.
An example of the case $n = 3$ and $m = 2$ is depicted in Figure \ref{fig : definition of cones}.
\begin{figure}[h]\label{fig : Coxeter-Dynkin of finite ADE}
\begin{tikzpicture}
  \draw[black!50,-] (0, -3.2) -- (0, 3.2) node[above]{Im};
  \draw[black!50,-] (-3.2, 0) -- (3.2, 0) node[right]{Re};
  
  \fill[gray,opacity=0.3](0,0)to(2,3)to(3,3)to(3,1);
  \fill[gray,opacity=0.3](0,0)to(-2,-3)to(-3,-3)to(-3,-1);
  \fill[gray,opacity=0.3](0,0)to(-0.5,3)to(-3,3)to(-3,2.7);
  \fill[gray,opacity=0.3](0,0)to(0.5,-3)to(3,-3)to(3,-2.7);
  \fill[gray,opacity=0.3](0,0)to(-2.5,3)to(-3,3)to(-3,1);
  \fill[gray,opacity=0.3](0,0)to(2.5,-3)to(3,-3)to(3,-1);
  
  \draw[thick] (-3,-1) -- (3,1);
  \draw[thick] (-2,-3) -- (2,3);
  \draw[thick] (-0.5,3) -- (0.5,-3);
  \draw[thick] (-3,2.7) -- (3,-2.7);
  \draw[thick] (-2.5,3) -- (2.5,-3);
  \draw[thick] (-3,1) -- (3,-1);
  
  \draw[black] (-1.1,2.2) node {$Z(A_1^+)$};
  \draw[black] (1.1,-2.2) node {$Z(A_1^-)$};
  \draw[black] (2,1.5) node {$Z(A_2^+)$};
  \draw[black] (-1.9,-1.5) node {$Z(A_2^-)$};
  \draw[black] (-1.9,1) node {$Z(A_3^-)$};
  \draw[black] (2,-1.1) node {$Z(A_3^+)$};
  
  \filldraw[white, draw=black, thick] (0,0) circle (0.08);
  \draw[black!50,-] (8, -3.2) -- (8, 3.2) node[above]{Im};
  \draw[black!50,-] (4.8, 0) -- (11.2, 0) node[right]{Re};

  \fill[gray,opacity=0.3](8,0)to(7.5,3)to(5,3)to(5,1);
  \fill[gray,opacity=0.3](8,0)to(8.5,-3)to(11,-3)to(11,-1);
  \fill[gray,opacity=0.3](8,0)to(10,3)to(11,3)to(11,1);
  \fill[gray,opacity=0.3](8,0)to(6,-3)to(5,-3)to(5,-1);

  \draw[thick] (5,-1) -- (11,1);
  \draw[thick] (6,-3) -- (10,3);
  \draw[thick] (7.5,3) -- (8.5,-3);
  \draw[thick] (5,1) -- (11,-1);
  \draw[thick] (5,-0.4) -- (11,0.4);
  
  \draw[black] (6.5,1.5) node {$C^+_1$};
  \draw[black] (9.5,-1.5) node {$C^-_1$};
  \draw[black] (10,1.5) node {$C^+_2$};
  \draw[black] (6.1,-1.5) node {$C^-_2$};
  
  \draw[->] ([shift={(8,0)}]0:1.6) arc[radius=2, start angle=0, end angle=14.5];
  \draw[->] ([shift={(8,0)}]0:2) arc[radius=1, start angle=0, end angle=15];

  \draw[black] (10,0) node[below] {$\theta$};
  \draw[black] (9.6,0.05) node[below] {$\theta'$};
  
  \filldraw[white, draw=black, thick] (8,0) circle (0.08);
\end{tikzpicture}
\caption{
An example of the case $n = 3$ and $m = 2$. 
Left: The images of cones $A_i^\pm$. 
Right: The cones $C^\pm_j$ and phases $\theta$ and $\theta'$.
}
\label{fig : definition of cones}
\end{figure}
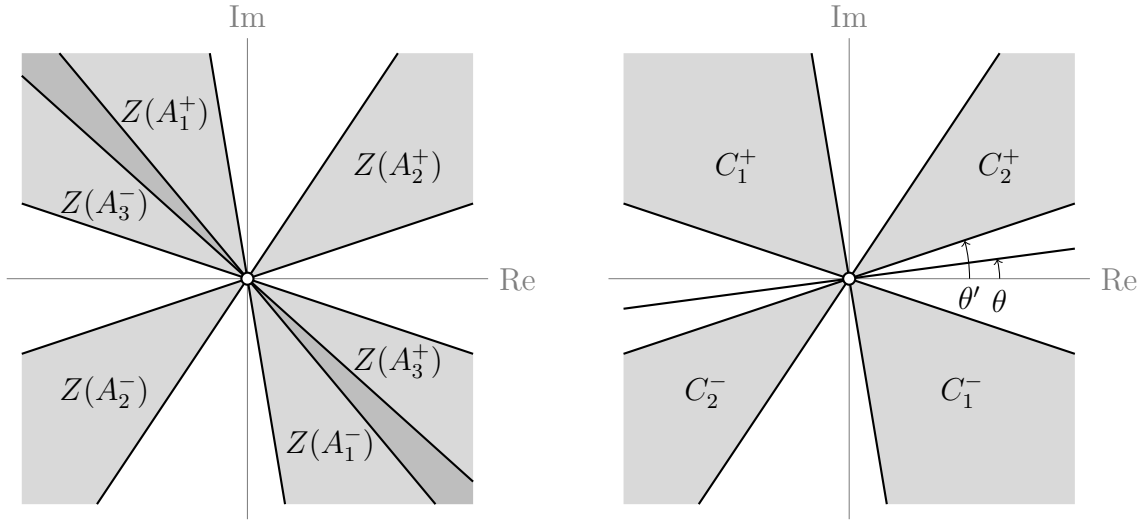
It follows from Lemma \ref{lem : intersections of cones} that $C^+_i \cap C^+_j = \emptyset$ for any $i \ne j$ and $C^+_i \cap C^-_j = \emptyset$ for any $i, j$.
Since $C^\pm_i$ is generated by $Z(\{ \alpha \in A^\pm_i \mid \| \alpha \| = 1\})$, it follows from the connectedness and compactness of $Z(\{ \alpha \in A^\pm_i \mid \| \alpha \| = 1\})$ that the set $C_i$ is a connected closed cone.
Hence, by Lemma \ref{lem : finiteness of connected components}, we have a decomposition of $Z(A(\sigma))$ by connected closed cones 
\begin{equation*}
Z(A (\sigma)) = C^+_1 \sqcup \dots \sqcup C^+_m \sqcup C^-_1 \sqcup \cdots \sqcup C^-_m .
\end{equation*}
Since each connected component is a connected closed cone, one can choose $\theta \in \RR$ so that $Z (A (\sigma)) \cap \RR e^{\pi \sqrt{-1} \theta}= \emptyset$.
By Lemma \ref{lem : connected closed cone in H}, each connected component $C^\pm_i$ included in $C^{(\theta, \theta + 1]}$ has the minimum phases $\phi^-_{C^\pm_i}$.
Put $\theta' \coloneqq \min \{ \phi^-_C \mid C = C^\pm_i, ~ C \subset C^{(\theta, \theta + 1]} \}$.
(As an example, see the right picture of Figure \ref{fig : definition of cones}.)
By construction, we have $\theta' > \theta$ and $Z(A(\sigma)) \cap C^{(\theta, \theta')} = \emptyset$.
\qed
\end{pf}
\begin{prop}\label{prop : existence of gaps}
There exist real numbers $\theta, \theta' \in \RR$ with $\theta' > \theta$ such that 
\begin{equation*}
Z(\C^\mathrm{ss} (\sigma) \cap \Delta) \cap C^{(\theta, \theta')} = \emptyset.
\end{equation*}
\end{prop}
\begin{pf}
By Lemma \ref{lem : fatten semistable imaginary cone} and Proposition \ref{prop : image of the fatten semistable imaginary cone}, we have real numbers $\theta$ and $\theta'$ such that 
\begin{equation*}
Z(\cl{\C^\mathrm{ss} (\sigma) \cap \Delta_\mathrm{im}}_0) \cap C^{(\theta, \theta')} = \emptyset.
\end{equation*}
By Lemma \ref{lem : limit of real rays}, the limit rays for the set $Z(\cl{\Delta_\mathrm{re}}_0)$ lie in $Z(\cl{\Delta_\mathrm{im}}_0)$.
Therefore, one can choose $\theta, \theta' \in \RR$ so that $Z(\cl{\C^\mathrm{ss} (\sigma) \cap \Delta}_0) \cap C^{(\theta, \theta')} = \emptyset$.
\qed
\end{pf}
\begin{rem}
Dimitrov--Haiden--Katzarkov--Kontsevich studied closed intervals $I \subset \RR$ with 
\begin{equation*}
Z(\C^\mathrm{ss} (\sigma) \cap \Delta) \cap C^{I} \ne \emptyset
\end{equation*}
from the viewpoint of the density of phases, as an analogue of the density of the set of slopes of closed geodesics on a Riemann surface.
For affine Dynkin quivers (resp., Dynkin quivers), Proposition \ref{prop : existence of gaps} also follows from \cite[Corollary 3.15]{DHKK} (resp., \cite[Lemma 3.13]{DHKK}).
\end{rem}
Let $\theta, \theta' \in \RR$ be as in Proposition \ref{prop : existence of gaps}. 
Since the derived category $\D^b(Q)$ is Krull--Schmidt, each $\sigma$-semistable object is decomposed into indecomposable $\sigma$-semistable objects.
It follows from Kac's Theorem (Corollary \ref{cor : Kac}) that for an indecomposable $\sigma$-semistable object $E$ we have $[E] \in \C^\mathrm{ss}(\sigma) \cap \Delta$. 
Therefore, Proposition \ref{prop : existence of gaps} yields $\P(\theta, \theta') = \{ 0 \}$.
By Proposition \ref{prop : criterion of algebraicity}, the heart $\P(\theta, \theta + 1]$ is algebraic.
We have finished the proof of Theorem \ref{thm : main 1}.
\qed
%

\subsection{Connectedness of the space of stability conditions}
In this section, we show the connectedness of the space of stability condition as a conclusion of Theorem \ref{thm : main 1}.
We first collect some notions and results concerning silting objects and simple-minded collections.
For more details, see \cite{AI, KY}.
Let $\D$ be a triangulated category.
For any objects $E, F \in \D$, we write $\Hom^p_\D(E, F) \coloneqq \Hom_\D(E, F[p])$ for $p \in \ZZ$ and $\Hom^\bullet_\D (E, F) \coloneqq \bigoplus_{p \in \ZZ} \Hom_\D (E_i, E_j [p])[-p]$.
An object $M \in \D$ is called {\em silting} if the following two conditions hold
\begin{itemize}
\item $\Hom^p_\D(M, M) = 0$ for all positive integers $p > 0$. 
\item The thick closure of $M$ is $\D$.
\end{itemize}
A silting object $M$ is called {\em tilting} if $\Hom^p_\D(M, M) = 0$ for $p \ne 0$.
We say two silting objects are {\em equivalent} if their additive closures coincide.
Let $M$ be a basic silting object in $\D$ and $M = M_1 \oplus \cdots \oplus M_\mu$ the decomposition into indecomposable objects.
For $i = 1, \dots, \mu$, the {\em left mutation} at the direct summand $M_i$ is the object $\mu_i^L (M) \coloneqq M'_i \oplus \bigoplus_{j \ne i} M_j$,
where $M'_i$ is the mapping cone of the minimal left approximation of $M_i$ with respect to $\bigoplus_{j \ne i} M_j$.
Similarly, one can define the {\em right mutation} $\mu_i^R (M)$ at the direct summand $M_i$ (see \cite{AI, KY} for more details).
Define a graph $\Silt(\D)$ as the oriented graph whose vertices are all equivalence classes of basic silting objects in $\D$ and whose edges correspond to left mutations.
A collection $X = \{ X_1, \dots, X_\mu \}$ of objects of $\D$ is said to be {\em simple-minded} if the following three conditions hold:
\begin{itemize}
\item $\Hom^p_\D (X_i, X_j) = 0$ for $p < 0$.
\item $\End_\D(X_i) \cong k$ and $\Hom_\D (X_i, X_j) = 0$.
\item The thick closure of $X$ is $\D$.
\end{itemize}
Note that a simple-minded collection is an unordered set.
One can define an equivalence relation of simple-minded collections in the natural way.
The {\em left mutation} $\mu_i^L (X)$ of $X = \{ X_1, \dots, X_\mu \}$ at $X_i$ is a new collection $\{ X'_1, \dots, X'_\mu \}$ such that $X_i' = X_i[1]$ and $X'_j$ for $j \ne i$ is the mapping cone of the left approximation of $X_j[-1]$ with respect to the closure of $X_i$.
Similarly, one can define the {\em right mutation} $\mu_i^R (X)$ at $X_i$ (see \cite{KY, KQ} for more details).
Define a graph $\SMC(\D)$ as the oriented graph whose vertices are all equivalence classes of simple-minded collections in $\D$ and whose edges correspond to left mutations.
Koenig--Yang established a remarkable correspondence among silting objects, simple-minded collections and algebraic hearts for finite-dimensional algebras.
Denote by $\per(Q)$ the perfect derived category of dg $k Q$-modules.
Note that we have $\per(Q) \cong \D^b(Q)$.
\begin{prop}[{\cite[Theorem 6.1 and 7.12]{KY}, cf.~\cite[Theorem 5.9]{KQ}}]\label{prop : correspondences}
Let $Q$ be a finite connected acyclic quiver.
\begin{enumerate}
\item There exists an isomorphism of oriented graphs between $\Silt(\per(Q))$ and $\EG^\mathrm{alg}(\D^b(Q))$:
\begin{equation*}
\Silt(\per(Q)) \overset{\cong}{\longrightarrow} \EG^\mathrm{alg}(\D^b(Q)), \quad M \mapsto \fmod(\End(M)).
\end{equation*}
\item There exists an isomorphism of oriented graphs between $\EG^\mathrm{alg}(\D^b(Q))$ and $\SMC(\D^b(Q))$:
\begin{equation*}
\EG^\mathrm{alg}(\D^b(Q)) \overset{\cong}{\longrightarrow} \SMC(\D^b(Q)), \quad \A \mapsto \Sim(\A),
\end{equation*}
\begin{equation*}
\SMC(\D^b(Q)) \overset{\cong}{\longrightarrow} \EG^\mathrm{alg}(\D^b(Q)), \quad X \mapsto \ex{X}.
\end{equation*}
\end{enumerate}
\qed
\end{prop}
\begin{rem}
For a finite-dimensional $k$-algebra $\Lambda$, Koenig--Yang also established one-to-one correspondences among bounded co-$t$-structures in $\per(\Lambda)$ and $\Silt(\per(\Lambda))$ and $\EG^\mathrm{alg}(\D^b(\Lambda))$.
\end{rem}
Aihara--Iyama studied the transitivity of silting mutations on several triangulated categories. 
As a conclusion of their result, we have 
\begin{prop}\label{prop : connectedness of the algebraic exchange graph}
Let $Q$ be a finite connected acyclic quiver.
The algebraic exchange graph $\EG^\mathrm{alg}(\D^b(Q))$ is connected.
\end{prop}
\begin{pf}
It was proved by \cite[Theorem 1.2]{AI} that $\Silt(\per(Q))$ is connected.
Hence, the statement follows from Proposition \ref{prop : correspondences}.
\qed
\end{pf}
As a consequence, every stability condition can be described by the $\CC$-action and iteration of simple tilts.
More precisely, we have the following:
\begin{thm}\label{thm : simple tilts and rotation}
Let $Q$ be a finite connected acyclic quiver.
For any stability condition on $\D^b(Q)$, there is a real number $\theta \in \RR$ such that the heart $\P(\theta, \theta + 1]$ is obtained from the standard heart $\fmod (k Q)$ by iteration of simple tilts.
In particular, the space $\Stab(\D^b(Q))$ is connected.
\end{thm}
\begin{pf}
It follows from Lemma \ref{lem : simple tilts} and Proposition \ref{prop : connectedness of the algebraic exchange graph} that $\Stab^\mathrm{alg}(\D^b(Q))$ is connected. 
Hence, the statement easily follows from Theorem \ref{thm : main 1}.
\qed
\end{pf}
It was proved by \cite{KV} (cf.~\cite[Appendix A]{Qiu1}) that for a Dynkin quiver $\vec{\Delta}$ any heart is obtained from the standard heart $\fmod (k \vec{\Delta})$ by iteration of simple tilts.
Especially, $\Stab(\D^b(\vec{\Delta}))$ is connected.
Theorem \ref{thm : simple tilts and rotation} can be regarded as a generalization of the result.
\begin{rem}
In \cite{CHQ}, the authors observed several sorts of connected components of spaces of stability conditions.
For a connected component $G$ of $\EG^\mathrm{alg}(\D)$ of a triangulated category $\D$, denote by $\Stab^\circ(\D)$ the connected component that contains subsets $U(\A)$ for $\A \in G$.
The component $\Stab^\circ(\D)$ is said to be of {\em finite type} if 
\begin{equation*}
\Stab^\circ(\D) = \bigcup_{\A \in G} U(\A).
\end{equation*}
$\Stab^\circ(\D)$ is said to be of {\em generic-finite type} if 
\begin{equation*}
\Stab^\circ(\D) = \CC \cdot \bigcup_{\A \in G} U(\A)
\end{equation*}
and it is not of finite type.
In their terminologies, Theorem \ref{thm : main 1} states that all components of $\Stab(\D^b(Q))$ are of generic-finite type.
Theorem \ref{thm : simple tilts and rotation} states that $\Stab(\D^b(Q))$ consists of a unique generic-finite type component.
\end{rem}
\bigskip
Next, we consider totally semistable stability conditions on $\D^b(Q)$.
One can construct a totally semistable stability condition with a hereditary algebraic heart (cf.~Section \ref{sec : gldim} and \cite[Lemma 5.1]{Qiu2}).
In order to construct another hereditary algebraic heart, we give a quick review of  Auslander--Platzeck--Reiten tiltings.
We refer to \cite{BB, APR} for more details.
Denote by $P_i$ and $S_i$ the indecomposable projective $k Q$-module and simple $k Q$-module corresponding to a vertex $i \in Q_0$, respectively.
Define the {\em Brenner--Butler tilting module} $T_i$ by 
\begin{equation*}
T_i \coloneqq \tau^{-1} S_i \oplus \bigoplus_{j \ne i} P_i,
\end{equation*}
where $\tau \in \Aut(\D^b(Q))$ is the Auslander--Reiten translation.
Similarly, one can define the (dual) Brenner--Butler tilting module $T^\vee_i$ with respect to the vertex $i \in Q_0$.
If the vertex $i \in Q_0$ is sink (resp., source), then $T_i$ (resp., $T^\vee_i$) is called the {\em Auslander--Platzeck--Reiten} tilting (APR-tilting for short).
An APR-tilting is interpreted by the {\em Bernstein--Gelfand--Ponomarev} reflection (BGP reflection for short) of the quiver $Q$.
More precisely, it was proved by \cite{APR} that 
\begin{equation*}
\End (T_i) \cong k (\mu^L_i (Q)), \quad \End (T^\vee_i) \cong k (\mu^R_i (Q)), \quad 
\end{equation*}
where $\mu^L_i$ (resp., $\mu_i^R$) is the BGP reflection with respect to the sink (resp., source) $i \in Q_0$.
Note that an APR-tilt is a simple tilt.
On the other hand, it is well-known that a tilting object induces a derived equivalence between $\D^b(Q) \cong \D^b(\End(T))$.
Hence, an acyclic quiver that is  obtained from $Q$ by iteration of BGP reflections is derived equivalent to the original acyclic quiver $Q$.
Moreover, the converse statement also holds.
Namely, the following result is known:
\begin{prop}[{\cite[Section 4.8]{Hap}}]\label{prop : GBP reflection and hereditary hearts}
Let $Q$ and $Q'$ be finite connected acyclic quivers.
We have $\D^b(Q) \cong \D^b(Q')$ if and only if $\fmod(k Q')$ is obtained from $\fmod(k Q)$ by iteration of APR-tilts.
\qed
\end{prop}
By the above proposition, one can deduce the following description of hearts arising from totally semistable stability conditions.
\begin{cor}\label{cor : Toss and APR tilt}
Let $Q$ be a finite connected acyclic quiver. 
Assume that a stability condition $\sigma = (Z, \P)$ on $\D^b(Q)$ is totally semistable.
Then, there exists $\theta \in \RR$ such that the heart $\P(\theta, \theta + 1]$ is obtained from the standard heart $\fmod(k Q)$ by iteration of APR-tilts.
Moreover, every algebraic heart $\A$ of the form $\P(\theta, \theta + 1]$ for some $\theta \in \RR$ arises in this way.
\end{cor}
\begin{pf}
Let $\A$ be an algebraic heart of the form $\P(\theta, \theta + 1]$ for some $\theta \in \RR$.
Note that Theorem~\ref{thm : main 1} ensures the existence of such an algebraic heart.
By Proposition \ref{prop : correspondences}, there exists a silting object $M \in \Silt (\per(Q))$ such that $\A \cong \fmod(\End(M))$.
It follows from Lemma \ref{lem : gldim for hearts} that $\gldim{\A} \le 1$, in particular the $k$-algebra $\End(M)$ is hereditary.
Hence, there exists a finite connected acyclic quiver $Q'$ such that $\End(M) \cong k Q'$ (e.g.,~see \cite{ASS}).
Note that since the rank of $K_0(\A)$ is $\mu$, the number of vertices of $Q'$ is also $\mu$.
The simple-minded collection $X = \{ X_1, \dots, X_\mu \}$ of $\D^b(Q)$ corresponding to $M$ is given as simple $k Q'$-modules.
Since $\per(Q) \cong \D^b(Q)$, it follows from \cite[Lemma 5.2]{KY} that the indecomposable projective $k Q'$-modules $P_1, \dots, P_\mu$ satisfies $M \cong P_1 \oplus \cdots \oplus P_\mu$, which yields $M$ is tilting.
Therefore, since we have $\D^b(Q) \cong \D^b(Q')$, the statement follows from Proposition \ref{prop : GBP reflection and hereditary hearts}.
\qed
\end{pf}
As a direction for further research, we are interested in the contractibility conjecture for spaces of stability conditions. 
This conjecture is related to the classical $K(\pi,1)$-conjecture, since certain hyperplane arrangements can be realized as quotients of the stability spaces of some Calabi--Yau categories (see \cite{Bri, Bri3}). 
Therefore, the contractibility conjecture may be viewed as a categorical analogue of the $K(\pi,1)$-conjecture. 
In \cite{QW}, the authors established the contractibility of stability spaces of $\D^b(\vec{\Delta})$ and $\D_{fd}(\Gamma_N \vec{\Delta})$ for Dynkin quivers, where $\Gamma_N \vec{\Delta}$ is the {\em $N$-Calabi--Yau completion} of $\vec{\Delta}$. 
Based on works \cite{Qiu2, Qiu3, QW}, we expect the following:
\begin{conj}
Let $Q$ be a finite connected acyclic quiver.
\begin{enumerate}
\item $\Stab(\D^b(Q))$ contracts to $\mathrm{Toss}(\D^b(Q))$.
\item $\mathrm{Toss}(\D^b(Q))$, the set of totally semistable stability conditions on $\D^b(Q)$, is contractible.
\end{enumerate}
In particular, $\Stab(\D^b(Q))$ is contractible.
\end{conj}
This conjecture was proved for the case of the affine $A_{p,q}$-quiver $Q = A^{(1)}_{p,q}$ (see \cite{HKK} and \cite{QZ}).
For affine Dynkin quivers, the contractibility of $\mathrm{Toss}(\D^b(Q))$ is also proved in \cite{QZ}.
We expect Corollary \ref{cor : Toss and APR tilt} is helpful to prove the contractibility of $\mathrm{Toss}(\D^b(Q))$.

\subsection{Full $\sigma$-exceptional collections}
We recall related notions of an exceptional collection.
Let $\D$ be a triangulated category.
\begin{itemize}
\item An object $E \in \D$ is called {\em exceptional} if $\Hom^\bullet_\D(E,E) \cong k$.
\item An ordered set $\E = (E_1, \dots, E_\mu)$ consisting of exceptional objects $E_1, \dots, E_\mu$ is called {\em exceptional collection} if $\Hom^p_\D(E_i, E_j) \cong 0$ for all $p \in \ZZ$ and $i > j$.
\item An exceptional collection $\E$ is called {\em full} if the smallest full triangulated subcategory of $\D$ containing all elements in $\E$ is equivalent to $\D$ as a triangulated category.
\item An exceptional collection $\E = (E_1, \dots, E_\mu)$ is called {\em Ext} if $\Hom^p_\D (E_i, E_j) \cong 0$ for $i \ne j$ and $p \le 0$.
\item An exceptional collection $\E = (E_1, \dots, E_\mu)$ is called {\em monochromatic} if for any $i, j = 1, \dots, \mu$, the $\ZZ$-graded $\CC$-vector space $\Hom^\bullet_\D (E_i, E_j) \not\cong 0$ is concentrated in a single degree.
\end{itemize}
For any full exceptional collection $(E_1, \dots, E_\mu)$, since $\D$ is of finite type one can choose integers $p_1, \dots, p_\mu \in \ZZ$ so that the shifted full exceptional collection $(E_1 [p_1], \dots, E_\mu [p_\mu])$ is Ext.
In \cite{Mac}, Macr\`{i} proved that the extension closure of a full Ext-exceptional collection forms an algebraic heart.
In the case of the derived category of an acyclic quiver, one can obtain the inverse statement.
\begin{prop}\label{prop : algebraic heart and full Ext exceptional collection}
Let $Q$ be a finite connected acyclic quiver.
For each algebraic heart $\A$ in $\D^b(Q)$, there is a monochromatic full Ext-exceptional collection $\E$ such that the extension closure of $\E$ is $\A$.
\end{prop}
\begin{pf}
By Proposition \ref{prop : correspondences}, we have a simple minded collection $\{ X_1, \dots, X_\mu \}$ such that $\A = \ex{X_1, \dots, X_\mu}$.
It was shown by \cite[Proposition 3.14]{IJ} that one can choose an order so that $(X_1, \dots, X_\mu)$ forms a full exceptional collection.
By the definition of a simple-minded collection, the full exceptional collection $(X_1, \dots, X_\mu)$ is Ext.
It was shown by \cite[Proposition 6.4]{KQ} that an algebraic heart obtained from the standard heart by iteration of simple tilts is monochromatic.
Hence, it follows from Proposition \ref{prop : connectedness of the algebraic exchange graph} that the full exceptional collection $(X_1, \dots, X_\mu)$ is monochromatic.
\qed
\end{pf}
By Proposition \ref{prop : algebraic heart}, one can consider stability conditions associated with a full Ext-exceptional collection.
Conversely, Dimitrov--Katzarkov introduced the notion of a full $\sigma$-exceptional collection with respect to a stability condition $\sigma$.
\begin{defn}[{\cite[Definition 3.17]{DK1}}]\label{defn : sigma-exceptional collection}
Let $\sigma = (Z, \P) \in \Stab(\D)$ be a stability condition on $\D$.
An exceptional collection $\E = (E_1, \dots, E_\mu)$ in $\D$ is called {\em $\sigma$-exceptional collection} if the following three properties hold:
\begin{itemize}
\item For each $i = 1, \dots, \mu$, the object $E_i$ is $\sigma$-semistable.
\item $\E$ is an Ext-exceptional collection.
\item There exists a real number $\theta \in \RR$ such that $\theta - 1 < \phi(E_i) \le \theta$ for $i = 1, \dots, \mu$.
\end{itemize}
\end{defn}
The existence of full $\sigma$-exceptional collections has been studied in several cases.
For generalized Kronecker quivers, Macr\`{i} showed that every stability condition admits a full $\sigma$-exceptional collection \cite[Lemma 4.2]{Mac} (cf.~\cite[Lemma A.1]{DK1}).
The same statement was proved for the affine $A_2$-quiver $A^{(1)}_{1,2}$ by \cite[Theorem 10.1]{DK1} and \cite[Theorem 4.16]{RW}. 
For Dynkin quivers, the same result was also shown by \cite[Theorem 1.2]{Ota}.
In \cite[Conjecture 3.11]{Ota}, it is conjectured that the same results hold for extended Dynkin quivers.
The following theorem not only gives an affirmative answer to this conjecture but also generalizes these known results.
\begin{thm}\label{thm : stability condition and full exceptional collections}
Let $Q$ be a finite acyclic quiver.
Every stability condition $\sigma = (Z, \P)$ on $\D^b(Q)$ admits a monochromatic full $\sigma$-exceptional collection.
\end{thm}
\begin{pf}
The statement follows from Theorem \ref{thm : main 1} and Proposition \ref{prop : algebraic heart and full Ext exceptional collection}.
\qed
\end{pf} 
Based on the existence of a full $\sigma$-exceptional collection, in \cite{DK2} and \cite{DK3}, Dimitrov--Katzarkov studied the contractibility of $\Stab(\D^b(Q))$ for the affine $A_2$-quiver and the generalized Kronecker quiver, respectively.
We believe that Theorem \ref{thm : stability condition and full exceptional collections} will play an important role to prove the contractibility of $\Stab(\D^b(Q))$.
Motivated by mirror symmetry, it is expected that the space of stability conditions on a triangulated category has a certain (natural) Frobenius structure in some settings.
From the viewpoint of singularity theory, it is natural to study the relation between stability condition on a derived directed Fukaya category and full exceptional collections.
Based on the correspondence between singularities and (generalized) root systems, we expect that Theorem \ref{thm : stability condition and full exceptional collections} will play an important role in constructing a conjectural Frobenius structure, which should be isomorphic to another one obtained by the invariant theory of the Weyl group, on the space of stability conditions on a derived category of an acyclic quiver.


\end{document}